\documentclass[12pt]{article}
\usepackage[cp1251]{inputenc}
\usepackage[english]{babel}
\usepackage {amsmath}
\usepackage {amsfonts}
\usepackage {amssymb}
\usepackage{graphicx}

\begin{document}   

\newtheorem{theorem}{Theorem}
\newtheorem{lemma}{Lemma}
\newtheorem{opr}{Defenition}
\centerline{\Large \bf On arithmetical nature of Tichy-Uitz's function}
\vspace*{4mm}
\centerline{\Large \bf Jabitskaya E.\,N. 
\footnote{ Research is supported by RFBR grant No 09-01-00371-a}}
\vspace*{4mm}

In \cite{TU}
R.\,F.\,Tichy and J.\,Uitz introduced a one parameter family 
$g_{\lambda}$,
$\lambda \in (0,1)$,
of singular  functions. 
When $\lambda=1/2$ the function $g_{\lambda}$ coincides with
the famous Minkowski question mark function. In this paper we describe
the arithmetical nature of the function $g_{\lambda}$ when 
$\lambda = \frac{3-\sqrt{5}}{2}$. 



\section{Stern-Brocot sequences}
We remind the definition of Stern-Brocot sequences $\mathcal{F}_n, 
n = 0,1,2,\ldots$. 

Consider the two-point set $\mathcal{F}_0 = \left\{ \frac 0 1 , \frac 1 1 \right\}$.
Let $n \geqslant 0$ and 
$$\mathcal{F}_n = 
\left\{0 = x_{0,n} < x_{1,n}< \ldots < x_{N(n),n} = 1 \right\},$$
where $x_{j,n} = p_{j,n}/q_{j,n}$, $(p_{j,n},q_{j,n}) = 1$, 
$j = 0,\ldots, N(n)$ and $N(n) = 2^n$. 
Then
$$\mathcal{F}_{n+1} =\mathcal{F}_{n} \cup Q_{n+1}$$ 
with 
$$Q_{n+1} = \left\{x_{j-1,n}\oplus x_{j,n},\quad j = 1, \ldots, N(n)\right\}.$$
Here
$$\frac{a}{b} \oplus \frac{c}{d} = \frac{a+b}{c+d}$$
is the mediant of the fractions $\frac{a}{b}$ and $\frac{c}{d}$.

The elements of $Q_n$ can be characterized in the following way.
A rational number $\xi \in [0,1]$ belongs to $Q_n$ if and only if in the 
continued fraction expansion
\begin{equation}\label{cf_real}
\xi = [ 0; a_1, a_2, \ldots, a_m] = 
0+\cfrac{1}{a_1+\cfrac{1}{a_{2} 
+ \ldots + \cfrac{1}{a_m}}},   \quad a_j \in \mathbb{N},\, a_m \geqslant 2.
\end{equation}
the sum of 
partial quotients is exactly $n+1$:
$$
S(\xi ) : = a_1+...+a_m = n+1.$$
So $\mathcal{F}_n$ consists of all rational $\xi \in [0,1]$ such that
$ S(\xi) \leqslant n+1$.

\section{Tichy-Uitz's singular functions}

In \cite{TU}
R.\,F.\,Tichy and J.\,Uitz considered a one parameter family 
$g_{\lambda}$,
$\lambda \in (0,1)$,
of singular  functions. 
In this section we describe the construction of $g_{\lambda}$ from \cite{TU}.
This construction is an inductive one.

Given $\lambda \in (0,1)$ put
$$
g_{\lambda} (0) =g_{\lambda} (0/1) = 0, \quad
g_{\lambda} (1) =g_{\lambda} (1/1) = 1.
$$ 
Suppose that $g_{\lambda} (x)$ is defined for all elements $x \in \mathcal{F}_n$.
Then we
define $g_{\lambda}(x) $ for $ x \in Q_{n+1}$.
Each $x\in Q_{n+1}$ is of the form
 $x=x_{j-1,n}\oplus x_{j,n}$ where $x_{j-1,n}$ and $x_{j,n}$ are consecutive 
elements from $\mathcal{F}_n$.
Then
$$
g_{\lambda} (x_{j-1,n}\oplus x_{j,n}) = g_{\lambda} (x_{j-1,n}) +
\left(g_{\lambda}(x_{j,n}) - g_{\lambda}(x_{j-1,n}) \right) \lambda.
$$
So we have defined $g_{\lambda}$ for all rational numbers from $[0,1]$.
One can  see that the  function $g_{\lambda}(x) $ is a
continuous function from $\mathbb Q \cap [0,1]$ to $[0,1]$.
So it can be extended to a continuous function from the whole segment $[0,1]$
to $[0,1]$.

For every $\lambda$ the function $g_{\lambda}(x) $ increases in $ x\in [0,1]$.
By the Lebesgue theorem $g_{\lambda}(x) $ is a differentiable function
almost everywhere. Moreower, it is easy to see that 
$g_{\lambda}'(x) =0$ almost everywhere (is the scence of Lebesgue mesure).
Certain properties of functions $g_{\lambda}(x) $ were investigated in
\cite{TU}. Some related topics one can find in \cite{DZ} and \cite{TG}.
Here we should note that in the case $\lambda = 1/2$ the function
$ g_{1/2}(x)$ coincides  with the famous Minkowski question mark function
$?(x)$.This function may be considered as the limit distribution function for  
Stern-Brocot  sequences $\mathcal{F}_n$.
The aim of the present paper is to  explain the arithmetical nature of the function
$g_{\lambda}(x) $ when $\lambda = \frac{3-\sqrt{5}}{2}$.

\section{Minkowski's function $?(x)$}
Let us consider the function $g_{1/2} (x) = ?(x)$.
This function was introduced by Minkowski. As it follows from the definition
of $g_{\lambda}$ for $\lambda=1/2$: 
$$
? (0) = ? (0/1) = 0, \quad
? (1) = ? (1/1) = 1.
$$
and for $x_{j-1,n}, x_{j,n} \in \mathcal{F}_n$
$$
? (x_{j-1,n}\oplus x_{j,n}) = \frac{? (x_{j-1,n}) + ?(x_{j,n})}{2}. 
$$
The definition of $?(x)$ for irrational $x$ follows by continuity.

R.\,Salem in \cite{Salem} found a new presentation for $?(x)$.
If $x \in (0,1)$ is represented in the 
form of regular continued fraction 
\begin{equation}\label{cfrac}
x = [ 0; a_1, a_2, \ldots, a_m, \ldots] = \cfrac{1}{a_1+\cfrac{1}{a_{2} 
+ \ldots + \cfrac{1}{a_m+ \frac{1}{\ldots}}}},
\end{equation}
then
\begin{equation}\label{formula_salem}
?(x) = \frac{1}{2^{a_1 - 1}} - \frac{1}{2^{a_1 + a_2 -1}} + 
\frac{1}{2^{a_1 + a_2 + a_3 - 1}}- \ldots
\end{equation}
For rational $x$ the representation (\ref{cfrac}) and consequently
(\ref{formula_salem}) is finite.

Minkowski's question mark function may be treated as the limit distribution
function for Stern-Brocot sequences in the following sense:
\begin{equation} \label{distr_func}
?(x) = \lim_{n \to \infty} \frac{\sharp \left\{\xi \in \mathcal{F}_n:
\xi \leqslant x\right\}}{\sharp \mathcal{F}_n} = 
\lim_{n \to \infty} \frac{\sharp \left\{\xi \in \mathcal{F}_n:
\xi \leqslant x\right\}}{2^n + 1}.
\end{equation}
A finite formula for the right hand side of (\ref{distr_func}) was
given by T. Rivoal in the preprint \cite{RI}.
Various properties of Minkowski question mark function were investigated in papers
\cite{Denjoy} by A.\,Denjoy, \cite{Para} by P.\,Viader, J.\,Paradis, L.\,Bibiloni
and in \cite{4} by A.\,A.\,Dushistova, I.\,D.\,Kan and N.\,G.\,Moshchevitin. 

\section{General form of formula (\ref{formula_salem})}
The formula (\ref{formula_salem}) can be generalized on the whole family of
functions $g_{\lambda}$ in the following way.

\noindent {\bf Proposition}
Let $x, \lambda \in (0,1)$ and 
$x = [0;a_1, \ldots, a_m, \ldots]$ is the regular continued fraction expansion
of $x$, then
\begin{multline}\label{general_formula_salem}
g_{\lambda} (x) = {\lambda}^{a_1 - 1} -{\lambda}^{a_1 - 1} (1-\lambda)^{a_2 } +
{\lambda}^{a_1-1} (1-\lambda)^{a_2 } \lambda^ {a_3}  - \ldots 
+\\+(-1)^{m+1} {\lambda}^{\sum\limits_{\left({1 \leqslant i \leqslant m}
,\, {i \equiv 1 \mod 2}\right)} a_i- 1}
(1-\lambda)^{\sum\limits_{\left({1 \leqslant i \leqslant m},\,
{i \equiv 0 \mod 2}\right)} a_i} + \ldots.
\end{multline}
{\bf Proof}:
By definition of $g_{\lambda}$ 
$$g_{\lambda} (0) = 0,\quad g_{\lambda} (1) = 1$$
and
\begin{equation}\label{f1}
g_{\lambda} (x_{j-1,n}\oplus x_{j,n}) = g_{\lambda} (x_{j-1,n}) +
\left(g_{\lambda}(x_{j,n}) - g_{\lambda}(x_{j-1,n}) \right) \lambda,
\end{equation}
where $x_{j-1,n}$ and $x_{j,n}$ are consecutive 
elements from $\mathcal{F}_n$.
We can also rewrite the formula (\ref{f1}) in the form
\begin{equation}\label{f2}
g_{\lambda} (x_{j-1,n}\oplus x_{j,n}) = g_{\lambda} (x_{j,n}) -
\left(g_{\lambda}(x_{j,n}) - g_{\lambda}(x_{j-1,n}) \right) (1-\lambda).
\end{equation}

The equality
$$g_{\lambda}(1/a_1) = \lambda^{a_1-1}$$
follows from the formula (\ref{f1}) immediately since
$1/a_1 = \underbrace{0 \oplus \ldots \oplus 0 \oplus }_{(a_1-1)\,\, \text{times}} 1$.
Suppose that the formula (\ref{general_formula_salem}) is proved for
$x = [0;a_1, \ldots, a_m]$, then it is enough to prove it for 
$y = [0;a_1, \ldots, a_m+1]$
and for $z = [0;a_1, \ldots, a_m, 2]$.

Let $m$ is odd, then by applying formula (\ref{f1}) we get
\begin{multline}
g_{\lambda}(y) = g_{\lambda} ([0;a_1, \ldots, a_{m-1}] \oplus
[0;a_1, \ldots, a_{m}])= \\ = 
g_{\lambda} ([0;a_1, \ldots, a_{m-1}]) + \lambda (g_{\lambda} 
([0;a_1, \ldots, a_{m}])-g_{\lambda} ([0;a_1, \ldots, a_{m-1}]))= \\ =
g_{\lambda} ([0;a_1, \ldots, a_{m-1}]) + 
{\lambda}^{\sum\limits_{\left({1 \leqslant i \leqslant m}
,\,{i \equiv 1 (\rm{mod} 2)}\right)} a_i- 1}
(1-\lambda)^{\sum\limits_{\left({1 \leqslant i \leqslant m}
,\,{i \equiv 1 \mod 2}\right)} a_i}
\lambda,
\end{multline}
and by applying formula (\ref{f2}) we get
\begin{multline}
g_{\lambda}(z) = g_{\lambda} ([0;a_1, \ldots, a_m+1] 
\oplus [0;a_1, \ldots, a_m]) = \\ =
g_{\lambda} ([0;a_1, \ldots, a_m]) - (1-\lambda)(g_{\lambda}([0;a_1, \ldots, a_m]) - 
g_{\lambda}([0;a_1, \ldots, a_m+1])= \\ =
g_{\lambda} ([0;a_1, \ldots, a_m]) -
{\lambda}^{\sum\limits_{\left({1 \leqslant i \leqslant m}
,\,{i \equiv 1 \mod 2}\right)} a_i- 1}
(1-\lambda)^{\sum\limits_{\left({1 \leqslant i \leqslant m}
,\,{i \equiv 1 \mod 2}\right)} a_i}
(1-\lambda)^2.
\end{multline}

For even $m$ the proof is analogously.

\section{Regular reduced continued fractions and the main result}
Any real number $x$ can be expressed uniquely in the 
form 
\begin{equation}\label{rrcfrac}
x = [[ b_0; b_1, b_2, \ldots, b_l, ...]] = b_0-\cfrac{1}{b_1-\cfrac{1}{b_{2} 
- \ldots - \cfrac{1}{b_l - \frac{1}{\ldots}}}}, \quad b_i \geqslant 2,
\end{equation}
which is known as regular reduced continued 
fraction (eine reduziert-regelmassige Kettenbrouch \cite{Finkelshtein}, 
\cite{Perron}).

For a rational number 
$x \in (0,1)$ the representation (\ref{rrcfrac}) takes the form:
\begin{equation} \label{rrcf_real}
x = [[1; b_1, \ldots, b_l]].
\end{equation}
For such $x$ we denote $L(x) = b_1 + \ldots + b_l$.

Analogously to the sequence $\mathcal{F}_n$ we define the sequence $\Xi_n$:
$$\Xi_n := \{0,1\} \cup \left( \mathop{\cup}\limits_{1 \leqslant k \leqslant n}
\Theta_k \right),$$
where 
$\Theta_k = \{ x \in Q : L(x) = k + 1\}$, $k \geqslant 1$.

We arrange the elements of $\Xi_n$ in the increasing order:
$$\Xi_k = \{0 = \xi_{1,n} < \xi_{2,n} < \ldots < \xi_{\sharp \Xi_n,n} = 1\}.$$

We would like to note that in the special case $\lambda = \tau^2$ formula 
(\ref{general_formula_salem}) gives:
\begin{multline}\label{formula2}
g_{\tau^2} (x) = {\tau}^{2a_1 - 2} -{\tau}^{2a_1 + a_2 - 2} +
{\tau}^{2a_1 + a_2 + 2a_3- 2}  - \ldots 
+\\+(-1)^{m+1} {\tau}^{\sum_{i = 1}^{m} \alpha_i a_i- 2} + \ldots,
\end{multline}
where 
\begin{equation*}
\alpha_m = 
\left\{
\begin{array}{ll}
1, &\text{if }m\text{ is even}, \\
2, &\text{if }m\text{ is odd}.  
\end{array}
\right.
\end{equation*} 
For rational $x$ the representation (\ref{formula2}) is finite.

The Theorem 1 below is the main results of the present paper.
It generalizes the formula (\ref{distr_func}) on the
regular reduced continued fractions.
\begin{theorem}\label{theorem1}
Function $g_{\lambda}$, where $\lambda = \tau^2 = \frac{3-\sqrt{5}}{2}$, 
$\tau = \frac{\sqrt{5}-1}{2}$ coincides with the distributional function
of the sequence $\Xi_n$, that is
$$
g_{\tau^2} (x) =\lim_{n \to \infty} \frac{\sharp \{\xi \in \Xi_n: 
\xi \leqslant x\}}{\sharp \Xi_n}, \quad x \in (0,1).$$
\end{theorem}

Now we consider the function
$${\mathcal{M}} (x) := \lim_{n \to \infty} \frac{\sharp \{\xi \in \Xi_n: 
\xi \leqslant x\}}{\sharp \Xi_n}, \quad x \in (0,1).$$
Our purpose is to prove that $\mathcal{M} (x) = g_{\lambda}$.
Function $\mathcal{M} (x)$ is increasing as a distribution function,
so it is enough to prove that $\mathcal{M}(x)$ coincides 
with $g_{\tau^2}(x)$ for rational $x$, that is
\begin{equation}\label{aa}
{\mathcal{M}} (x\oplus y) = {\mathcal{M}} (x) + \left({\mathcal{M}} (y) 
- {\mathcal{M}} (x)\right) \tau^2.
\end{equation}
for any two consecutive elements of $\Xi_n$ for any $n$.

\section{Auxiliary results}
The following result is well known. We present it without a proof.
\begin{lemma}\label{lemma0}
Let $x$ is represented in the form (\ref{cf_real})
and in the form (\ref{rrcf_real}). To get the set $(b_1, \ldots, b_l)$
from $(a_1, \ldots, a_m)$ we should replace $a_i$ by
\begin{enumerate}
\item $\underbrace{2\ldots2}_{a_i - 1}$ if $i$ is odd
(empty string if $a_i = 1$).
\item $a_i + 2$ if $i$ is even and $i \ne m$.
\item $a_i + 1$ if $i$ is even and $i = m$.
\end{enumerate}
\end{lemma}
\begin{lemma}
For the number of elements in $\Theta_n$ one has
$$\sharp \Theta_1 = 1,\,
\sharp \Theta_2 = 1,\,
\sharp \Theta_{n+1} = \sharp \Theta_{n}+ \sharp \Theta_{n-1},$$ 
so 
$\sharp \Theta_n = F_{n},$ where 
$$F_n = \frac{\left(\frac{1+\sqrt{5}}{2}\right)^n - 
\left(\frac{1-\sqrt{5}}{2}\right)^n}{\sqrt{5}}$$
is the nth Fibonacci number.
\end{lemma}
{\bf Proof.}
We prove the lemma by induction.
Since $\Theta_1 = \{1/2\}$, $\Theta_2 = \{2/3\}$, then
the base of induction is true.
Let us suppose that the lemma is true for $k \leqslant n$ and
$x=[[1; b_1, \ldots, b_l]] \in \Theta_{n+1}$, then
$b_1 + \ldots + b_l = n + 2$. There are two cases: either $b_l = 2$ or 
$b_l \geqslant 2$. In the first case $b_1 + \ldots+ b_{l-1} = n$, so
$[[1; b_1, \ldots, b_{l-1}]] \in \Theta_{n-1}$, in the second case
$b_1 + \ldots + b_{l}-1 = n+1$, so $[[1; b_1, \ldots, b_{l}-1]] \in \Theta_n$.
Thus we have one-one correspondence between $\Theta_{n-1} \cup \Theta_{n}$
and $\Theta_{n+1}$, and so $\sharp \Theta_{n+1} = \sharp \Theta_{n}+ \sharp \Theta_{n-1}$.
\begin{opr}
Let $x, y, z$ be consecutive elements of $\Xi_n$, $y \in \Theta_n$.
We denote the mediant $x \oplus y$ by
$y^l$, the mediant $y\oplus z$ we denote
by $y^r$.
\end{opr}
\begin{lemma}
Let $x, y, z$ be consecutive elements of $\Xi_n$, $y \in \Theta_n$,
then $y^l \in \Theta_{n+2}$, $y^r \in \Theta_{n+1}$. 
\end{lemma}
{\bf Proof.}
Let $y = [[1; b_1, \ldots, b_{s}]]$. Then $y^l = [[1; b_1, \ldots, b_{s},2]]$,
$y^r = [[1; b_1, \ldots, b_{s}+1]]$.

Now let us construct an infinite tree $D$ whose
nodes are labeled by rationals in $(0,1)$. We identify the nodes 
with the rationals they labeled by.
The root is labeled by $1/2$. From node $x$ come two arrows: 
the left arrow goes to $x^l$ and the right arrow goes to $x^r$.
The nodes of the tree $D$ are partitioned into levels. 
$1/2$ belongs to the level 1.
If $x$ belongs to the level $n$, then $x^r$ belongs to the level $n+1$,
and $x^l$ belongs to the level $n+2$ (figure 1).
\begin{figure}[h]\label{tree}
\centerline{\includegraphics{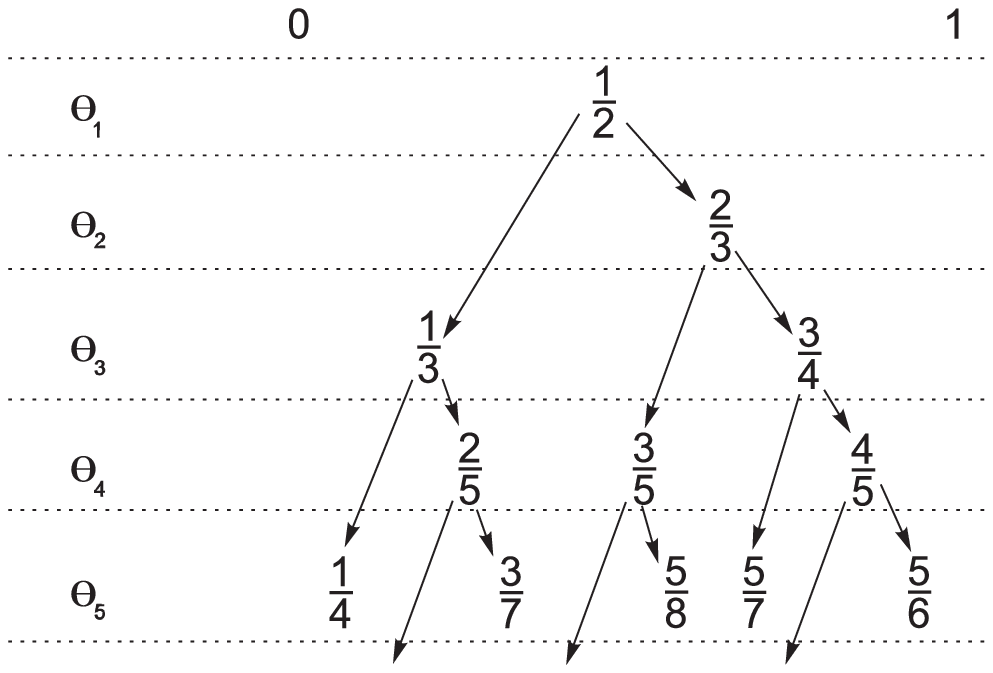}}
\caption{}
\end{figure}

It follows from the construction of the tree that nodes from level $n$
of $D$ are marked by numbers from $\Theta_n$.
So $x$ belongs to the level $n$ if and only if $x \in \Theta_n$.

The subtree of $D$ with root in the node $x$ we denote by $D^{(x)}$.
The set of nodes of $D$ from level 1 to level $n$ we denote by $D_n$.
The set of nodes of $D^{(x)} \cap D_n$ we denote by $D^{(x)}_{n}$.
Note that there exist a levels preserving isomorphism between $D$ and $D^{(x)}$. 
If $x$ belongs to the level $n$, then
$$
\sharp D^{(x)}_{m} = \sharp D_{m-n+1}.
$$
Besides
$$
\sharp D_n =\sharp \Theta_1 +\sharp \Theta_2 + \ldots +\sharp \Theta_n = 
F_1 + F_2 + \ldots + F_n = F_{n+2} - 1.
$$ 

\section{Proof of Theorem 1}
We remind that it is enough to prove (\ref{aa}) for any consecutive 
elements of $\Xi_n$ $x$ and $y$. 

To prove the equality (\ref{aa}) we consider the subtree $D^{(x \oplus y)}$
of $D$. Note that
$$\left\{\xi \in D^{(x \oplus y)}\right\}\cup \{y\} 
= \{\xi \in Q: x < \xi \leqslant y\}.$$
Consequently 
$$
{\mathcal{M}} (y) - {\mathcal{M}} (x) = \lim_{m \to \infty}
\frac{\sharp \{\xi \in \Xi_m: x < \xi \leqslant y\}}
{\sharp \Xi_m} = 
\lim_{m \to \infty} \frac{\sharp D^{(x \oplus y)}_m}
{\sharp D_m}.
$$
On the other hand 
$$
{\mathcal{M}} (x \oplus y) - {\mathcal{M}} (x) = \lim_{m \to \infty}
\frac{\sharp \{\xi \in \Xi_m: x < \xi \leqslant x \oplus y\}}
{\sharp \Xi_m} = 
\lim_{m \to \infty} 
\frac{\sharp {D}^{(x \oplus y)^l}_m}{\sharp D_m}.
$$

Let $x \oplus y \in \Theta_{k}$, then $(x \oplus y)^l \in \Theta_{k+2}$.
Therefore 
\begin{multline*}
\frac{{\mathcal{M}} (x \oplus y) - {\mathcal{M}} (x)}{{\mathcal{M}} (y) - 
{\mathcal{M}} (x)} = \lim_{m \to \infty} \frac{\sharp D^{(x \oplus y)^l}_m}
{\sharp D^{(x \oplus y)}_m} = \\=
\lim_{m \to \infty} \frac{\sharp D_{m-k-1}}{\sharp D_{m-k+1}} =
\lim_{m \to \infty} \frac{F_{m-k+1}}{F_{m-k+3}} = \tau^2.
\end{multline*}

\end{document}